\documentclass[12pt,twoside]{amsart}
\usepackage{amsfonts,amsmath}
\usepackage[all,knot,arc]{xy}

\def\endpf{\relax\ifmmode\expandafter\endproofmath\else
  \unskip\nobreak\hfil\penalty50\hskip.75em\hbox{}\nobreak\hfil\bull
  {\parfillskip=0pt \finalhyphendemerits=0 \bigbreak}\fi}
\def\bull{\vbox{\hrule\hbox{\vrule\kern3pt\vbox{\kern6pt}\kern3pt\vrule}\hrule}}

\newtheorem{defn}{Definition}[section]
\newtheorem{lemma}[defn]{Lemma}
\newtheorem{sublemma}[defn]{Sublemma}
\newtheorem{theorem}[defn]{Theorem}

\newtheorem{remark}[defn]{Remark}
\newtheorem{proposition}[defn]{Proposition}
\newtheorem{corollary}[defn]{Corollary}

\newtheorem{maintheorem}{Theorem}
\newtheorem{mainproposition}[maintheorem]{Proposition}

\newcommand{\zz}{{\mathbb Z}}
\newcommand{\rr}{{\mathbb R}}
\newcommand{\nn}{{\mathbb N}}
\newcommand{\qq}{{\mathbb Q}}

\renewcommand{\aa}{{\mathbb A}}
\newcommand{\bb}{{\mathbb B}}

\newcommand{\ii}{{\mathbf i}}
\newcommand{\jj}{{\mathbf j}}
\newcommand{\kk}{{\mathbf k}}

\newcommand{\froyshov}{Fr{\o}yshov}
\newcommand{\ozsvath}{Ozsv\'{a}th}
\newcommand{\szabo}{Szab\'{o}}
\newcommand{\spin}{\ifmmode{\rm Spin}\else{${\rm spin}$\ }\fi}
\newcommand{\spinc}{\ifmmode{{\rm Spin}^c}\else{${\rm spin}^c$\ }\fi}
\newcommand{\spinct}{\mathfrak t}
\newcommand{\spincs}{\mathfrak s}
\newcommand{\tors}{{\it Tors}}

\newcommand{\calt}{\mathcal{T}}

\newcommand{\detQ}{\delta}

\newcommand{\dt}{{\,.\,}}
\newcommand{\ds}{\displaystyle}
\def\ceil#1{\left\lceil#1\right\rceil}
\def\floor#1{\left\lfloor#1\right\rfloor}

\newcommand{\HFhat}{\widehat{HF}}

\newenvironment{narrow}[2]{%
 \begin{list}{}{%
  \setlength{\topsep}{0pt}%
  \setlength{\leftmargin}{#1}%
  \setlength{\rightmargin}{#2}%
  \setlength{\listparindent}{\parindent}%
  \setlength{\itemindent}{\parindent}%
  \setlength{\parsep}{\parskip}%
 }%
\item[]}{\end{list}}

\newif\ifpic
\picfalse    
\pictrue   

\DeclareMathOperator\im{Im}
\DeclareMathOperator\Hom{Hom}

\begin{document}

\title{A characterisation of the $\zz^n\oplus\zz(\delta)$ lattice 
and definite nonunimodular intersection forms}
\author{Brendan Owens and Sa\v{s}o Strle}
\date{\today}
\thanks{B. Owens was supported in part by NSF grant DMS-0604876 and by the LSU Council on Research Summer Stipend Program.\\
S.Strle was  supported in part by the MSZS of the
Republic of Slovenia research program No. P1-0292-0101-04 and research
project No. J1-6128-0101-04.}

\begin{abstract}
We prove a generalisation of Elkies' theorem to nonunimodular definite 
forms (and lattices). Combined with inequalities of \froyshov~ and of \ozsvath\ 
and \szabo, this gives a simple test of whether a rational homology 3-sphere 
may bound a definite four-manifold.  As an example we show that small positive
surgeries on torus knots do not bound negative-definite
four-manifolds.
\end{abstract}

\maketitle

\pagestyle{myheadings}
\markboth{BRENDAN OWENS AND SA\v{S}O STRLE}{A CHARACTERISATION OF THE
 $\zz^n\oplus\zz(\detQ)$ LATTICE}


\section{Introduction}
\label{sec:intro}

The intersection pairing of a smooth compact four-manifold, possibly with boundary, is
an integral symmetric bilinear form $Q_X$ on $H_2(X;\zz)/\tors$; 
it is nondegenerate if the boundary 
of $X$ has first Betti number zero.  Characteristic covectors for $Q_X$ are elements of
$H^2(X;\zz)/\tors$ represented by the first Chern classes of \spinc structures.  In the case that
$Q_X$ is positive-definite there are inequalities due to \froyshov~ (using Seiberg-Witten
theory, see \cite{f}) and \ozsvath~ and \szabo~  
which give lower bounds on the
square of a characteristic covector.  It is helpful in this context to prove existence of
characteristic covectors with small square.  The following is the main result of this paper.

\begin{maintheorem}
\label{thm:mainthm}
Let $Q$ be an integral positive-definite symmetric bilinear form of rank $n$ and 
determinant $\detQ$. 
Then there exists a characteristic covector $\xi$ for $Q$ with 
$$\xi^2 \le 
\left\{\begin{array}{ll}
n-1 + 1/\detQ &  \mbox{if $\detQ$ is odd,}\\
n-1 &  \mbox{if $\detQ$  is even;}
\end{array}\right. 
$$ 
moreover the inequality is strict unless
$Q=(n-1)\langle1\rangle\oplus\langle\detQ\rangle$.
\end{maintheorem}

For unimodular forms this theorem was proved by Elkies \cite{elkies}. To prove 
Theorem \ref{thm:mainthm} we reinterpret the statement in terms of integral
lattices. In Section \ref{sec:lattice} we introduce the necessary notions and then
study short characteristic covectors of some special types of lattices. 
The proof of the theorem is presented in Section \ref{sec:proof}. The main idea 
is to use the theory of linking pairings to embed four copies of the lattice 
into a unimodular lattice and then apply Elkies' theorem along with results 
of Section \ref{sec:lattice}.

For unimodular forms, a well-known and useful 
constraint on the square of a characteristic vector is that $\xi^2$ is congruent
modulo 8 to the signature of the form.
In Section \ref{sec:congruence} we give a generalisation of this to nonunimodular forms.

In Section \ref{sec:app} we combine Theorem \ref{thm:mainthm}
with an inequality of \ozsvath~ and \szabo~ \cite[Theorem 9.6]{os4} to obtain the following theorem,
where $d(Y,\spinct)$ denotes the correction term invariant defined in \cite{os4}
for a \spinc structure $\spinct$ on a rational homology three-sphere $Y$:

\begin{maintheorem}
\label{thm:YboundsX}
Let $Y$ be a rational homology sphere with $|H_1(Y;\zz)|=\detQ$.
If Y bounds a negative-definite four-manifold $X$, and if either $\detQ$ is 
square-free or there is no torsion in $H_1(X;\zz)$, then
$$\max_{\spinct\in\spinc(Y)}4d(Y,\spinct)\ge
\left\{\begin{array}{ll}
1-1/\detQ & \mbox{if $\detQ$ is odd,}\\
1 & \mbox{if $\detQ$ is even.}
\end{array}\right.$$
The inequality is strict unless the intersection form of $X$ is 
$(n-1)\langle-1\rangle\oplus\langle -\detQ\rangle$.
\end{maintheorem}


As an application we consider manifolds obtained as integer surgeries on
knots.  Knowing which positive surgeries on a knot bound negative-definite manifolds has implications for existence of fillable contact structures and for the unknotting number of the knot  (see e.g.\ \cite{ls,bsinprep}).  
In particular we consider torus knots; it is well known that $(pq-1)$-surgery on 
the torus knot $T_{p,q}$ is a lens space. 
It follows that all large ($\ge pq-1$) surgeries on torus knots bound
both positive- and negative-definite four-manifolds (see e.g.\ \cite{bsinprep}).
We show
in Section \ref{sec:torus-knots} that small positive surgeries on torus
knots cannot bound negative-definite four-manifolds. In particular, we obtain the following
result; for a more precise statement see Corollary \ref{cor:torus-neg-def}.

\begin{mainproposition}
\label{prop:torus-neg-def}
Let $2\le p<q$ and $1\le n\le (p-1)(q-1)+2$. Then $+n$-surgery on $T_{p,q}$
cannot bound a negative-definite four-manifold with no torsion in $H_1$. 
\end{mainproposition}

\vskip2mm
\noindent{\bf Acknowledgements.}  
It is a pleasure to thank Noam Elkies who 
introduced to us the notion of
gluing of lattices, which is used in the proof of Proposition \ref{prop:embed};
John Morgan, who demonstrated a version of Corollary \ref{cor:mod4} to the first author;
Peter \ozsvath, who explained to us a proof
of Theorem \ref{thm:L-space-knot-d}; and  Ed Swartz, for a helpful conversation regarding the ``money-changing problem".  


We started this project  while we were Britton Postdoctoral Fellows at McMaster University, and some of the work was carried out while the first author was visiting Cornell University.


\section{Short characteristic covectors in special cases}
\label{sec:lattice}

An {\it integral lattice} $L$ of rank $n$ is the free abelian group $\zz^n$ along with
a nondegenerate symmetric bilinear form $Q:\zz^n\times\zz^n\to\zz$.
Let $V$ denote the tensor product of $L$ with $\rr$; thus $V$ is the vector space
$\rr^n$ to which the form $Q$ extends, and $L$ is a sublattice of $V$.  
The {\it signature} $\sigma(L)$ of the lattice $L$ is the signature of $Q$.
We say $L$ is {\it definite} if $|\sigma(L)|=n$.
For convenience denote $Q(x,y)$ by $x\cdot y$ and $Q(y,y)$ by $y^2$.
A set of generators $v_1,\ldots,v_n$ for $L$
forms a basis of $V$ satisfying $v_i\cdot v_j\in\zz$. With respect to such a basis
the form $Q$ is represented by the matrix $[v_i\cdot v_j]_{i,j=1}^n$. The 
{\it determinant} (or {\it discriminant}) of $L$
is the determinant of the form (or the corresponding matrix).

The {\it dual lattice} $L'\cong\Hom(L,\zz)$ consists of all vectors $x\in V$ satisfying
$x\cdot y\in\zz$ for all $y \in L$.  A {\it characteristic covector} for $L$
is an element $\xi\in L'$ with $\xi\cdot y\equiv y^2\pmod2$ for all $y\in L$.

We say a lattice is {\it complex} if it admits an automorphism $\ii$ 
with $\ii^2$ given by multiplication by $-1$. 
A lattice is {\it quaternionic} if it admits an action of the 
quaternionic group $\{\pm1,\pm \ii,\pm \jj,\pm \kk\}$, with $-1$ acting 
by multiplication.  Note that the rank of a complex lattice is even, and the rank of 
a quaternionic lattice is divisible by 4.  For any lattice $L$, let $L^m$
denote the direct sum of $m$ copies of $L$.  There is a standard way
to make $L\oplus L$ into a complex lattice and $L^4$ into a quaternionic lattice; for
example, the quaternionic structure is given by
\begin{eqnarray*}
\ii:(x,y,z,w)&\mapsto&(-y,x,-w,z)\\
\jj:(x,y,z,w)&\mapsto&(-z,w,x,-y).
\end{eqnarray*}

Let $\zz(\detQ)$ denote the rank one lattice 
with determinant $\detQ$; in particular $\zz=\zz(1)$.


\begin{lemma}
\label{lem:selfglue}
Let $\detQ\in\nn$, and let $p$ be a prime. Let $L$ be an index $p$ sublattice of 
$\zz^{n-1}\oplus\zz(\detQ)$.
Then $L$ has a characteristic covector $\xi$ with 
$\xi^2 < n-1$ 
unless $L\cong\zz^{n-1}\oplus\zz(p^2\detQ)$.
\end{lemma}

\proof
We may assume $L$ contains none of the summands of $\zz^{n-1}\oplus\zz(\detQ)$; any such
summand of $L$ contributes 1 to the right-hand side of the inequality and
at most 1 to the left-hand side.
 
If $n=1$ then clearly $L\cong\zz(p^2\detQ)$.  
Now suppose $n>1$.  Let $\{e,e_1,\ldots,e_{n-2},f\}$ 
be a basis for $\zz^{n-1}\oplus\zz(\detQ)$, where $e,e_1,\ldots,e_{n-2}$ have square 
1 and $f^2=\detQ$.  Then multiples of $e$ give coset representatives of $L$ in 
$\zz^{n-1}\oplus\zz(\detQ)$; it follows that a basis for $L$ is given by
$\{pe,e_1+s_1 e,\ldots,e_{n-2}+s_{n-2}e,f+te\}$. Here $s_i,t$ are nonzero residues 
modulo $p$ in $[1-p,p-1]$, whose parities we may choose if $p$ is odd.
With respect to this basis, the bilinear form on $L$ has matrix
$$Q=\left(\begin{matrix}
p^2 & ps_1 & ps_2 & ps_3 & \dots & ps_{n-2} & pt\\
ps_1 & 1+s_1^2 & s_1s_2 & s_1s_3 & \dots & s_1s_{n-2} & s_1t\\
 & & & & & & \\
\vdots & \vdots & \vdots & \vdots & & \vdots & \vdots\\
 & & & & & & \\
ps_{n-2} & s_1s_{n-2} &s_2s_{n-2} & s_3s_{n-2} & \dots & 1+s_{n-2}^2 & s_{n-2}t \\
pt & s_1t & s_2t & s_3t & \dots & s_{n-2}t & \detQ+t^2\\
\end{matrix}\right),$$
with inverse
$$Q^{-1}=\left(\begin{matrix}
\frac{1+\sum s_i^2 +t^2/\detQ}{p^2} & -\frac{s_1}p & -\frac{s_2}p & -\frac{s_3}p & \dots 
& -\frac{s_{n-2}}p & -\frac{t}{p\detQ} \\
-\frac{s_1}p & 1 & 0 & 0 & \dots & 0 & 0\\
 & & & & & & \\
\vdots & \vdots & \vdots & \vdots & & \vdots & \vdots\\
 & & & & & & \\
-\frac{s_{n-2}}p & 0 & 0 & 0 & \dots & 1 & 0\\
-\frac{t}{p\detQ} & 0 & 0 & 0 & \dots & 0 & \frac1\detQ\\
\end{matrix}\right).$$
With respect to the dual basis, an element of the dual lattice $L'$ is represented 
by an $n$-tuple $\xi\in\zz^n$.  An $n$-tuple corresponds to a characteristic 
covector if its components have the same parity as the corresponding 
diagonal elements of $Q$.

Suppose now that $p$ is odd.  Choose $s_i$ to be odd for all $i$ and $t\equiv\detQ\pmod2$.
Then $\xi=(1,0,\dots,0)$ is a characteristic covector whose square satisfies
$$\xi^2=\frac{1+\sum s_i^2 +t^2/\detQ}{p^2}<n-1,$$
noting that $|s_i|,\ |t|\le p-1$.

Finally if $p=2$ then $\xi=(0,0,\dots,0,(1+\detQ) \bmod2)$ is a characteristic vector with
$\xi^2<n-1$.
\endproof

\begin{lemma}
\label{lem:pglue}
Let $\detQ\in\nn$ be odd, and let $p$ be a prime with $p=2$ or $p\equiv 1\pmod4$.
Let $M$ be an index $p$ complex sublattice of 
$\zz^{2n-2}\oplus\zz(\detQ)^2$.
Then $M$ has a characteristic covector $\xi$ with 
$\xi^2 < 2n-2$ 
unless $M\cong\zz^{2n-2}\oplus\zz(p\detQ)^2$. 
\end{lemma}

\proof
We may assume $M$ contains none of the summands of $\zz^{2n-2}\oplus\zz(\detQ)^2$.
Suppose first that $n=1$.  Let $\{e,\ii e\}$ be a basis for $\zz(\detQ)^2$ with $e^2=\detQ$.
Then $\{pe, \ii e+se\}$ is a basis for $M$, for some $s$.
Now
\begin{eqnarray*}
\ii(\ii e+se)=-e+s\ii e &\in& M\\
\implies -e-s^2e&\in&M,
\end{eqnarray*}
from which it follows that $p$ divides $1+s^2$.  The bilinear form on $L$ has
matrix
$$Q=\left(\begin{matrix}p^2\detQ & ps\detQ\\ ps\detQ & \detQ(1+s^2)\end{matrix}\right)
=p\detQ\left(\begin{matrix}p & s\\ s & \frac{1+s^2}p\end{matrix}\right)
\cong p\detQ I,$$
since any positive-definite unimodular form of
rank 2 is diagonalisable.  Thus in this case $M\cong \zz(p\detQ)^2$.

Now suppose $n>1$.  Let $\{e,e_1,\dots,e_{2n-3},f_1,f_2\}$ be a basis for 
$\zz^{2n-2}\oplus\zz(\detQ)^2$, where $e,e_i$ have square 1 and $f_i^2=\detQ$.
A basis for $M$ is given by $\{pe,e_i+s_i e,f_i+t_i e\}$, where $s_i,t_i$ are nonzero
residues modulo $p$ which we may choose to be odd integers in $[1-p,p-1]$. 
The matrix in this basis is
$$Q=\left(\begin{matrix}
p^2 & ps_1 & ps_2 & ps_3 & \dots & ps_{2n-3} & pt_1 & pt_2\\
ps_1 & 1+s_1^2 & s_1s_2 & s_1s_3 & \dots & s_1s_{2n-3} & s_1t_1 & s_1t_2\\
ps_2 & s_1s_2 & 1+s_2^2 & s_2s_3 & \dots & s_2s_{2n-3} & s_2t_1 & s_2t_2\\
 & & & & & & &\\
\vdots & \vdots & \vdots & \vdots & & \vdots & \vdots & \vdots\\
 & & & & & & \\
ps_{2n-3} & s_1s_{2n-3} & s_2s_{2n-3} & s_3s_{2n-3} & \dots & 1+s_{2n-3}^2 & s_{2n-3}t_1 & 
s_{2n-3}t_2 \\
pt_1 & s_1t_1 & s_2t_1 & s_3t_1 & \dots & s_{2n-3}t & \detQ+t_1^2 & t_1t_2\\
pt_2 & s_1t_2 & s_2t_2 & s_3t_2 & \dots & s_{2n-3}t & t_1t_2 & \detQ+t_2^2\\
\end{matrix}\right),$$
with inverse
$$Q^{-1}=\left(\begin{matrix}
\frac{1+\sum s_i^2 +\sum t_i^2/\detQ}{p^2} & -\frac{s_1}p & -\frac{s_2}p & \dots 
& -\frac{s_{2n-3}}p & -\frac{t_1}{p\detQ} & -\frac{t_2}{p\detQ}\\
-\frac{s_1}p & 1 & 0 & \dots & 0 & 0 & 0\\
-\frac{s_2}p & 0 & 1 & 0 & \dots & 0 & 0\\
 & & & & & &\\
\vdots & \vdots & \vdots & & & \vdots\\
 & & & & & &\\
-\frac{s_{2n-3}}p & 0 & 0 & \dots & 1 & 0 & 0\\
-\frac{t_1}{p\detQ} & 0 & 0 & \dots & 0 & \frac1\detQ & 0\\
-\frac{t_2}{p\detQ} & 0 & 0 & \dots & 0 & 0 & \frac1\detQ\\
\end{matrix}\right).
$$

If $p=2$ then $Q$ is even, and $\xi=0$ is characteristic.

Suppose now that $p$ is odd. 
Then $\xi=(k,l_1,\dots,l_{2n-1})$ is characteristic
if $k$ is odd and each $l_i$ is even.  By completing squares we have that
$$\xi^2=\frac{1}{p^2}\left(k^2+\sum_{i=1}^{2n-3}(ks_i-pl_i)^2
+\frac1\detQ\sum_{i=1}^{2}(kt_i-pl_{2n-3+i})^2\right).$$
It clearly suffices to prove the statement for $\detQ=1$ and $n=2$; thus the following 
sublemma completes the proof.

\begin{sublemma}
Let $\ds F(k,l_1,l_2,l_3)=k^2+\sum_{i=1}^3(ks_i-l_ip)^2$,
where $s_i$ are odd numbers with $|s_i|\le p-2$.  Then there exist $k$ odd and $l_i$ even
for which $F<2p^2$.
\end{sublemma}
\proof
Let $K=\{2-p,4-p,\dots,p-2\}$, and let 
$$K_i=\left\{k\in K\ \bigg|\ \exists\ l_i\ \mbox{even with}\ |ks_i-l_ip|<\frac{p-1}{2}\right\}.$$
Note that $|K_i|=\frac{p-1}2$ (since $p\equiv1\pmod4$).

If there exists $k\in\ds\cap_{i=1}^3 K_i$ then with this choice of $k$ we find
$$F< p^2+3\left(\frac{p-1}2\right)^2<2p^2.$$
Otherwise let $K_{ij}$ denote $K_i\cap K_j$.
Then 
$$p-1\ge|\cup_{i=1}^3 K_i|=\sum_{i=1}^3 |K_i|-\sum_{i<j} |K_{ij}|,$$
from which it follows that $\sum_{i<j} |K_{ij}|\ge\frac{p-1}2$.  Thus there
exists some $k\in K_{ij}$ with $|k|\le \frac{p+1}2$.  With this $k$ we find
$$F< p^2 + 2\left(\frac{p-1}2\right)^2+\left(\frac{p+1}2\right)^2<2p^2.$$
\endproof

\begin{lemma}
\label{lem:qglue}
Let $\detQ\in\nn$ be odd, and let $q$ be a prime with $q\equiv3\pmod4$.
Let $N$ be an index $q^2$ quaternionic sublattice of 
$\zz^{4n-4}\oplus\zz(\detQ)^4$, with the quotient group $(\zz^{4n-4}\oplus\zz(\detQ)^4)/N$ having exponent $q$.
Then $N$ has a characteristic covector $\xi$ with 
$\xi^2 < 4n-4$ 
unless $N\cong\zz^{4n-4}\oplus\zz(q\detQ)^4$.
\end{lemma}

\proof
We may assume $N$ contains none of the summands of $\zz^{4n-4}\oplus\zz(\detQ)^4$.
Suppose first that $n=1$.  Let $\{e,\ii e,\jj e,\kk e\}$ be a basis for 
$\zz(\detQ)^4$ with $e^2=\detQ$. Note that $e$ and $\ii e$ represent generators 
of the quotient $\zz/q\oplus\zz/q;$ otherwise we have $e+s\ii e\in N$ for some $s$, 
and multiplication by $\ii$ yields $s^2+1\equiv 0\pmod q$, but $-1$ is not a 
quadratic residue modulo $q$.  Now a basis for $N$ is given by 
$\{qe,q\ii e,\jj e+s_1e+t_1\ii e,\kk e+s_2e+t_2\ii e\}$.  The quaternionic symmetry yields
$$s_1\equiv t_2,\ \ s_2\equiv -t_1,\ \ 1+s_1^2+s_2^2\equiv0\quad \pmod q;$$
it follows that the bilinear form $Q$ on $N$  factors 
as $q\detQ$ times a unimodular form.  Thus $N\cong \zz(q\detQ)^4$.

We now assume $n=2$ and $\detQ=1$ (the proof for any other $n,\detQ$ follows  
from this case).
Let $\{e,\ii e,v_1=\jj e,v_2=\kk e,v_3=f,v_4=\ii f,v_5=\jj f,v_6=\kk f\}$ 
be a basis of unit vectors for $\zz^8$.  Then
$\{qe,q\ii e,v_1+s_1e+t_1\ii e,\dots,v_6+s_6e+t_6\ii e\}$
is a basis for $N$.  In this case the quaternionic symmetry yields
$$ s_{i+1}\equiv -t_{i},\ \ t_{i+1}\equiv s_{i},\ \ 
1+s_1^2+s_2^2\equiv 0\quad \pmod q$$
for $i=1,3,5$, and at most one of $s_3,\dots,s_6$ is divisible by $q$.
We may choose $s_i,t_i$ in $[1-q,q-1]$ with $s_i+t_i$ odd.
Computation of $Q,Q^{-1}$ in the above basis for $N$ now yields that
$\xi=(k_1,k_2,l_1,\dots,l_6)$ is characteristic if $k_1,k_2$ are odd and $l_i$ are even, and
$$\xi^2=\frac1{q^2}\left(k_1^2+k_2^2+\sum_{i=1}^6(k_1s_i+k_2t_i-l_iq)^2\right).$$
There are now two cases to consider, depending on whether or not one of the $s_i$ is
zero.  Existence of a characteristic $\xi$ with $\xi^2<4$ follows in each case
from one of the following sublemmas.

\begin{sublemma}
Let F=$\ds k_1^2+k_2^2+\sum_{i=1}^6(k_1s_i+k_2t_i-l_iq)^2,$ where $s_i$ are odd and $t_i$ are
even integers in $[1-q,q-1]$.  Then there exist $k_1,k_2$ odd and $l_i$ even for which $F<4q^2$.
\end{sublemma}
\proof
Set $k_2=1$.
Let $K=\{2-q,4-q,\dots,q-2\}$, and let 
$$K_i=\left\{k\in K\ \bigg|\ \exists\ l_i\ \mbox{even with}\ 
|ks_i+t_i-l_iq|\le\frac{q-1}{2}\right\}.$$
Note that $|K_i|\ge\frac{q-1}2$ (in fact $K_i$ either contains $\frac{q-1}2$ or $\frac{q+1}2$
elements, depending on the value of $t_i$).

If there is a $k$ which is in the intersection of 4 of the $K_i$'s, then setting $k_1=k$ we
obtain
$$F\le 1+3q^2+4\left(\frac{q-1}2\right)^2<4q^2.$$
Otherwise let $K_{ijk}$ denote the triple intersection of $K_i,K_j,K_k$, and let
$K_{ij}'$ denote $K_i\cap K_j-\cup_k K_{ijk}$. Then
$$q-1\ge|\cup_{i=1}^6 K_i|=\sum_{i=1}^6 |K_i|-\sum_{i<j} |K_{ij}'|-2\sum_{i<j<k}|K_{ijk}|,$$
from which it follows that 
$$\sum_{i<j} |K_{ij}'|+2\sum_{i<j<k}|K_{ijk}|\ge2(q-1),$$ and hence
$$\sum_{i<j<k}|K_{ijk}|\ge\frac{q-1}2.$$
Thus there exists some $k\in K_{ijk}$ with $|k|\le\frac{q-1}2$.  Setting $k_1=k$ again yields
$$F\le 1+3q^2+4\left(\frac{q-1}2\right)^2<4q^2.$$\endproof

\begin{sublemma}
Let $F$=$\ds k_1^2+k_2^2+\sum_{i=1}^5(k_1s_i+k_2t_i-l_iq)^2+(k_2t_6-l_6q)^2,$ where
$s_i$ and $t_6$ are odd, and $t_1,\dots,t_5$ are even integers in $[1-q,q-1]$.  
Then there exist $k_1,k_2$ odd and $l_i$ even for which $F<4q^2$.
\end{sublemma}
\proof
First choose $k_2$ (and $l_6$) with $|k_2|,|k_2t_6-l_6q|\le \frac{q-1}2$.
Again let $K=\{2-q,4-q,\dots,q-2\}$, and for each $i\le5$ let
$$K_i=\left\{k\in K\ \bigg|\ \exists\ l_i\ \mbox{even with}\ 
|ks_i+k_2t_i-l_iq|\le\frac{q-1}{2}\right\};$$
then $|K_i|\ge\frac{q-1}2$.

If there is a $k$ which is in the intersection of 4 of the $K_i$'s, then setting $k_1=k$ we
obtain
$$F\le 2q^2+6\left(\frac{q-1}2\right)^2<4q^2.$$
Otherwise (with notation as above) we have
$$q-1\ge|\cup_{i=1}^5 K_i|=\sum_{i=1}^5 |K_i|-\sum_{i<j} |K_{ij}'|-2\sum_{i<j<k}|K_{ijk}|,$$
from which it follows that 
$$\sum_{i<j} |K_{ij}'|+2\sum_{i<j<k}|K_{ijk}|\ge\frac32(q-1),$$ and hence
$$\sum_{i<j<k}|K_{ijk}|\ge\frac{q+1}4.$$
Thus there exists some $k\in K_{ijk}$ with $|k|\le\frac{3(q+1)}4$.  Setting $k_1=k$ yields
$$F\le 2q^2+5\left(\frac{q-1}2\right)^2+\left(\frac{3(q+1)}4\right)^2<4q^2.$$\endproof


\section{Proof of the main theorem}
\label{sec:proof}

The bilinear form $Q$  on $L$ induces a symmetric bilinear $\qq/\zz$-valued pairing 
on the finite group $L'/L$, called the {\it linking pairing} associated to $L$.  
Such pairings $\lambda$ on a finite group $G$ were studied by Wall \cite{wall}
(see also \cite{kk}).  He observed that $\lambda$ splits into an orthogonal sum of pairings 
$\lambda_p$ on the $p$-subgroups $G_p$ of $G$. 
For any prime $p$ let $A_{p^k}$ (resp.\ $B_{p^k}$) denote the 
pairings on $\zz/p^k$ with $p^k$ times the square of a generator a quadratic
residue (resp.\ nonresidue) modulo $p$. 
If $p$ is odd then $\lambda_p$ 
decomposes into an orthogonal direct sum of these two types of pairings on cyclic subgroups.
The pairing on $G_2$ may be decomposed into cyclic summands (there are four equivalence
classes of pairings on $\zz/2^k$ for $k\ge 3$), plus two types of pairings on
$\zz/2^k\oplus\zz/2^k$; these are denoted $E_{2^k}, F_{2^k}$.  For an appropriate
choice of generators, in $E_{2^k}$ each cyclic generator has square 0, while in 
$F_{2^k}$ they each have square $2^{1-k}$.

\begin{proposition} 
\label{prop:embed}
Let $L$ be a 
lattice of rank $n$.
Then $L^4$ 
may be embedded as a quaternionic sublattice of
a unimodular quaternionic lattice $U$ of rank $4n$.
\end{proposition}
\proof
Consider an orthogonal decomposition of the linking pairing on $L'/L$
as described above.  Let $x\in L'$ represent a generator of a
cyclic summand $\zz/p^k$ with $k>1$.  Then $v=p^{k-1}x$ has $v^2\in\zz$
and $pv\in L$.  Adjoining $v$ to $L$ yields a lattice $L_1$ which contains
$L$ as an index $p$ sublattice (so that $\det L=p^2\det L_1$).
Similarly, if $x\in L'$ represents a generator of a cyclic summand
of $E_{2^k}$ or $F_{2^k}$ then $v=2^{k-1}x$ may be adjoined.
Finally if $x_1,x_2$ represent generators of two $\zz/2$ summands then
$x_1+x_2$ may be adjoined.
In this way we get
a sequence of embeddings
$$L=L_0\subset L_1 \subset \dots \subset L_{m_1},$$
where each $L_i$ is an index $p$ sublattice of $L_{i+1}$ for some prime $p$.
Moreover the linking pairing of $L_{m_1}$ decomposes into cyclic summands of prime order,
and the determinant of $L_{m_1}$ is either odd or twice an odd number.

Now let $M_0=L_{m_1}\oplus L_{m_1}$.  Note that $M_0$ is a complex sublattice
of $L'\oplus L'$.
Let $p$ be a prime which is either $2$ or is congruent to 1 modulo 4, so that
there exists an integer $a$ with $a^2\equiv-1\pmod p$.
Let $x\in L'$ represent a generator of a $\zz/p$ summand of the linking
pairing of $L_{m_1}$.  Then $v=(x,ax)\in L'\oplus L'$ has $v^2\in\zz$
and $pv\in M_0$.  Adjoining $v$ to $M_0$ yields a lattice $M_1$; 
since $\ii v+av\in M_0$, $M_1$ is preserved by $\ii$.  Continuing in this way 
we obtain a sequence of embeddings
$$M_0\subset M_1 \subset \dots \subset M_{m_2},$$
where each $M_i$ is an index $p$ sublattice of $M_{i+1}$ for some prime $p$
with $p=2$ or $p\equiv 1\pmod4$. Each $M_i$ is a complex sublattice of
$L'\oplus L'$.  We may arrange that each $M_i$ with $i>0$ has odd determinant.
The resulting linking pairing of $M_{m_2}$ is the orthogonal sum of pairings on cyclic
groups of prime order congruent to 3 modulo 4, and all the summands appear
twice.

Finally let $N_0=M_{m_2}\oplus M_{m_2}$.  This is a quaternionic sublattice of  
 $(L')^4$. Let $q\equiv3\pmod4$ be a prime, and suppose that $x\in L'$ 
generates a $\zz/q$ summand of the linking pairing of $M_{m_2}$.
%
There exist integers $a,b$ with $a^2+b^2\equiv-1\pmod q$ (take $m$ to be the smallest
positive quadratic nonresidue and choose $a,b$ so that $a^2\equiv m-1,\ b^2\equiv -m$).
Let $v_1=(x,0,ax,bx)$ and let $v_2=\ii v_1=(0,x,-bx,ax)$.  
Then $v_i^2, v_1\dt v_2\in\zz$, and $qv_i\in N_0$.  Adjoining $v_1,v_2$ 
to $N_0$ yields a quaternionic sublattice $N_1$ of $(L')^4$; note that 
$\jj v_1 + av_1-bv_2\in N_0$.  We thus obtain a sequence
$$N_0\subset N_1 \subset \dots \subset N_{m_3}=U,$$
where each $N_i$ is an index $q^2$ sublattice of $N_{i+1}$ for
some prime $q\equiv 3\pmod4$, $N_{i+1}/N_i$ has exponent $q$, each $N_i$ is 
quaternionic with odd determinant, and $U$ is unimodular.
\endproof

\begin{remark}
Using  the methods of Proposition \ref{prop:embed} 
one may show that $L^2$ embeds
in a unimodular lattice if and only if the prime factors of $\det L$ 
congruent to 3 mod 4 have even exponents.
Applying this to rank 1 lattices one recovers the classical fact that
any integer is expressible as the sum of 4 squares, and is the sum of 2 squares 
if and only if its prime factors congruent to 3 mod 4 have even exponents.
\end{remark}

\noindent{\it Proof of Theorem \ref{thm:mainthm}.}
Let $L$ be the lattice with form $Q$.
Embed $L^4$ in a unimodular lattice $U$ of rank $4n$ as in Proposition \ref{prop:embed}.
By Elkies' theorem \cite{elkies} either $U\cong\zz^{4n}$ or $U$ has a characteristic
covector $\xi$ with $\xi^2\le 4n-8$.  In the latter case, restricting to one of
the four copies of $L$ yields a covector $\xi|_L$ with $(\xi|_L)^2\le n-2$, which easily
satisfies the statement of the theorem.

This leaves the case that $U\cong\zz^{4n}$.  Let 
$L_0,\dots,L_{m_1},M_0,\dots,M_{m_2},N_0,\dots,N_{m_3}$ be as in the
proof of Proposition \ref{prop:embed}.
By successive application of Lemma \ref{lem:qglue}, we see that either 
$N_0\cong\zz^{4n-4}\oplus\zz(\detQ_1)^4$ for some $\detQ_1\in\nn$, 
or $N_0$ has a characteristic covector $\xi$
with $\xi^2<4n-4$.  In the latter case, restricting $\xi$ to one of the four copies of
$L$ yields $(\xi|_L)^2<n-1$.

If $N_0\cong\zz^{4n-4}\oplus\zz(\detQ_1)^4$ then 
$M_{m_2}\cong\zz^{2n-2}\oplus\zz(\detQ_1)^2$.  By Lemma \ref{lem:pglue} we
find that either $M_0\cong\zz^{2n-2}\oplus\zz(\detQ_2)^2$ or $M_0$ has a 
characteristic covector $\xi$
with $\xi^2<2n-2$.  In the latter case, restricting $\xi$ to one of the two copies of
$L$ embedded in $M_0$ yields $(\xi|_L)^2<n-1$.

Finally if $M_0\cong\zz^{2n-2}\oplus\zz(\detQ_2)^2$ then 
$L_{m_1}\cong\zz^{n-1}\oplus\zz(\detQ_2)$.  Successive application of
Lemma \ref{lem:selfglue} now yields that either
$L\cong\zz^{n-1}\oplus\zz(\detQ)$ or $L$ has a characteristic covector
$\xi$ with $\xi^2<n-1$.\endproof


\section{A congruence condition on characteristic covectors}
\label{sec:congruence}

Given a positive-definite symmetric bilinear form $Q$ 
of rank $n$ and determinant $\detQ$ with
$Q\ne (n-1)\langle1\rangle\oplus\langle\detQ\rangle$, one may ask
for an optimal upper bound on the square of a shortest characteristic covector
$\xi$.  The square of a characteristic covector of a unimodular form
is congruent to the signature modulo 8 (see for example \cite{serre}).
Thus if $\detQ=1$ we have $\xi^2\le n-8$.
In this section we give congruence conditions on the square of
characteristic covectors of forms of arbitrary determinant.
Lattices in this section are not assumed to be definite.
(The results in this section
may be known to experts but we have not found them in the literature.)

If $\xi_1,\xi_2$ are characteristic covectors of a
lattice $L$ of determinant $\detQ$ then their difference is
divisible by 2 in $L'$; it follows that
$\xi_1^2\equiv\xi_2^2$ modulo $\frac4\detQ$.
For a lattice with odd determinant the congruence holds modulo $\frac8\detQ$.
We will give a formula for this congruence class in terms of
the signature and linking pairing of $L$.  For a lattice with 
determinant $\detQ$  even or odd we
will determine the congruence class of $\xi^2$ modulo $\frac4\detQ$
in terms of the signature and determinant.

\begin{lemma}
\label{lem:evensig}
Suppose $\detQ=\ds\prod_{i=1}^r p_i^{k_i}\cdot\prod_{j=1}^s q_j^{l_j}$ 
where $p_i, q_j$ are odd primes, not necessarily
distinct.  Let $M$ be an even lattice 
of determinant $\detQ$ with linking pairing isomorphic to
$\ds\bigoplus_{i=1}^r A_{p_i^{k_i}}\oplus \bigoplus_{j=1}^s B_{q_j^{l_j}}.$
Then the signature of $M$ satisfies the congruence
$$\sigma(M)\equiv\ds\sum_{k_i\equiv1\,(2)}(1-p_i)+\sum_{l_j\equiv1\,(2)}(5-q_j)
\quad\mod8.$$
(For the definition of $A_{p^k}$ and 
$B_{p^k}$ see the beginning of Section \ref{sec:proof}.)
\end{lemma}
\proof
Let $G(M)$ denote the Gauss sum
$$\ds\frac1{\sqrt\detQ}\sum_{u\in M'/M}e^{i\pi u^2}.$$
(See \cite{taylor} for more details on Gauss sums and the Milgram Gauss sum formula.)
Then $G(M)$ depends only on the linking pairing of $M$ and in fact factors as
$$G(M)=\prod G(A_{p_i^{k_i}})\cdot\prod G(B_{q_j^{l_j}}).$$
The factors are computed in \cite[Theorem 3.9]{taylor} to be:
\begin{eqnarray*}
G(A_{p^k}) &=&
\left\{\begin{array}{ll}
1 &  \mbox{if $k$ is even,}\\
e^{2\pi i(1-p)/8} &  \mbox{if $k$  is odd;}
\end{array}\right.\\ 
G(B_{q^l}) &=&
\left\{\begin{array}{ll}
1 &  \mbox{if $l$ is even,}\\
e^{2\pi i(5-q)/8} &  \mbox{if $l$  is odd.}
\end{array}\right.\\ 
\end{eqnarray*}
(Note that in the notation of \cite{taylor},
$A_{p^k}$ corresponds to $(C_p(k);2)$ and
$B_{q^l}$ corresponds to $(C_q(l);2n_q)$,
where $n_q$ is a quadratic nonresidue modulo $q$.)

The congruence on the signature of $M$ now follows from the Milgram formula:
$$G(M)=e^{2\pi i\sigma(M)/8}.$$\endproof

\begin{proposition}
\label{prop:mod8}
Let $L$ be a lattice of 
odd determinant $\detQ$ with linking pairing isomorphic to 
$\ds\bigoplus_{i=1}^r A_{p_i^{k_i}}\oplus \bigoplus_{j=1}^s B_{q_j^{l_j}}.$
Let $\xi$ be a characteristic covector for $L$.  Then
$$\ds\xi^2\equiv \sigma(L)-\ds\sum_{k_i\equiv1\,(2)}(1-p_i)-\sum_{l_j\equiv1\,(2)}(5-q_j)
\quad\mod{\frac8\detQ}.$$
\end{proposition}
\proof
Let $M$ be an even lattice 
with the same linking pairing
as $L$.  (Existence of $M$ is proved by Wall in \cite{wall}.)  Then by \cite[Satz 3]{kp}
$L$ and $M$ are stably equivalent; that is to say, there exist unimodular lattices
$U_1, U_2$ such that $L\oplus U_1\cong M\oplus U_2$.  Let $\xi$ be a characteristic
covector for $L\oplus U_1$.  Then we have decompositions
$$\xi=\xi_L+\xi_{U_1}=\xi_M+\xi_{U_2}.$$
Taking squares we find
\begin{eqnarray*}
\xi_L{}\!^2&\equiv&\xi_M{}\!^2+\xi_{U_2}{}\!^2-\xi_{U_1}{}\!^2\\
&\equiv&\sigma(U_2)-\sigma(U_1)\\
&\equiv&\sigma(L)-\sigma(M)\\
&\equiv&\sigma(L)-\ds\sum_{k_i\equiv1\,(2)}(1-p_i)-\sum_{l_j\equiv1\,(2)}(5-q_j)
\quad\mod{\frac8\detQ},
\end{eqnarray*}
where the last line follows from Lemma \ref{lem:evensig}.
\endproof

\begin{corollary}
\label{cor:mod4}
Let $L$ be a lattice of 
determinant $\detQ\in\nn$, and
let $\xi$ be a characteristic covector for $L$.  Then
$$\xi^2\equiv
\left\{\begin{array}{ll}
\sigma(L)-1 + 1/\detQ &  \mbox{if $\detQ$ is odd}\\
\sigma(L)-1 &  \mbox{if $\detQ$  is even}
\end{array}\right.
\qquad\mod{\frac4\detQ}.
$$ 
\end{corollary}
\proof
If $\detQ$ is odd then Proposition \ref{prop:mod8} shows that
the congruence class of $\xi^2$ modulo $\frac4\detQ$ depends
only on the signature and determinant; the formula then follows
by taking $L$ to be the lattice with the form $r\langle1\rangle\oplus s\langle-1\rangle\oplus
\langle\detQ\rangle$ where $r+s=n-1$ and $r-s=\sigma(L)-1$.

If $\detQ$ is even then as in the proof of Proposition \ref{prop:embed} we find
that either $L$ or $L\oplus\zz(2)$ embeds as an index $2^k$ sublattice of
a lattice with odd determinant.  
It again follows that the congruence class of $\xi^2$ modulo $\frac4\detQ$ depends
only on the signature and determinant.
\endproof


\section{Proof of Theorem \ref{thm:YboundsX}}
\label{sec:app}

We begin by noting the following restatement of Theorem \ref{thm:mainthm}
for negative-definite forms:

\begin{theorem}
\label{thm:negdef}
Let $Q$ be an integral negative-definite symmetric bilinear form of rank $n$ and 
determinant of absolute value $\detQ$. 
Then there exists a characteristic covector $\xi$ for $Q$ with 
$$\xi^2 \ge 
\left\{\begin{array}{ll}
-n+1 - 1/\detQ &  \mbox{if $\detQ$ is odd,}\\
-n+1 &  \mbox{if $\detQ$  is even;}
\end{array}\right. 
$$ 
moreover the inequality is strict unless
$Q=(n-1)\langle-1\rangle\oplus\langle-\detQ\rangle$.
\end{theorem}

Let $Y$ be a rational homology three-sphere and $X$ a smooth
negative-definite four-manifold bounded by  $Y$, with $b_2(X)=n$.  For any $\spinc$
structure $\spinct$ on $Y$ let $d(Y,\spinct)$ denote the
correction term invariant of \ozsvath~ and \szabo~ \cite{os4}. It is
shown in \cite[Theorem 9.6]{os4} that for each $\spinc$ structure
$\spincs\in\spinc(X)$,
\begin{equation}
\label{eqn:thm9.6} 
c_1(\spincs)^2+n\le4d(Y,\spincs|_Y).
\end{equation}

The image of $c_1(\spincs)$ in $H^2(X;\zz)/\tors$ is a characteristic covector
for the intersection pairing $Q_X$ on $H_2(X;\zz)/\tors$.  
Let $\detQ$ denote the absolute value of the determinant of
$Q_X$ and $\im(\spinc(X))$ the image of the restriction map from $\spinc(X)$ to 
$\spinc(Y)$.  Then combining (\ref{eqn:thm9.6}) with Theorem \ref{thm:negdef} yields
$$\max_{\spinct\in\im(\spinc(X))}4d(Y,\spinct)\ge
\left\{\begin{array}{ll}
1-1/\detQ & \mbox{if $\detQ$ is odd,}\\
1 & \mbox{if $\detQ$ is even,}
\end{array}\right.$$
with strict inequality unless the intersection form of $X$ is 
$(n-1)\langle-1\rangle\oplus\langle-\detQ\rangle$.

Theorem \ref{thm:YboundsX} follows immediately since if either $\detQ$ is square-free
or if there is
no torsion in $H_1(X;\zz)$, then the restriction map from
$\spinc(X)$ to $\spinc(Y)$ is surjective, and $|H_1(Y;\zz)|=\detQ$.  

Similar reasoning yields the following variant of Theorem \ref{thm:YboundsX}:

\begin{proposition}
Let $Y$ be a rational homology sphere with $|H_1(Y;\zz)|=rs^2$, with $r$ square-free.
If Y bounds a negative-definite four-manifold $X$, then
$$\max_{\spinct\in\spinc(Y)}4d(Y,\spinct)\ge
\left\{\begin{array}{ll}
1-1/r & \mbox{if $r$ is odd,}\\
1 & \mbox{if $r$ is even.}
\end{array}\right.$$
\end{proposition}

\begin{remark}\label{thm:spinstr}
Suppose $Y$ and $X$ are as in the statement of Theorem \ref{thm:YboundsX} with $\detQ$ even.
If in fact $$\max_{\spinct\in\im(\spinc(X))}4d(Y,\spinct)=1,$$
then it is not difficult to see that this maximum must be attained at a spin structure.
\end{remark}


\section{Surgeries on $L$-space knots}
\label{sec:L-knots}

Let $K$ be a knot in the three-sphere and 
$$\Delta_K(T)=a_0+\sum_{j>0}a_j(T^j+T^{-j})$$ 
its Alexander polynomial. Torsion coefficients of $K$ are defined by 
$$t_i(K)=\sum_{j>0}ja_{|i|+j}\, .$$
Note that $t_i(K)=0$ for $i\ge N_K$, where $N_K$ denotes the degree of the
Alexander polynomial of $K$. For any $n\in\zz$ denote the $n$-surgery on $K$
by $K_n$. $K$ is called an {\em $L$-space knot} if for some $n>0$, $K_n$ is an $L$-space. 
Recall from \cite{lenssurg} that a rational homology sphere is called an {\em $L$-space} 
if $\HFhat(Y,\spincs)\cong\zz$ for each \spinc structure $\spincs$
(so its Heegaard Floer homology groups resemble those of a lens space).

There are a number of conditions coming from Heegaard Floer homology that an 
$L$-space knot $K$ has to satisfy. In particular (see \cite[Theorem 1.2]{lenssurg}), 
its Alexander polynomial has the form
\begin{equation}\label{eqn:L-knot-alex}
\Delta_K(T)=(-1)^k+\sum_{j=1}^k(-1)^{k-j}(T^{n_j}+T^{-n_j}).
\end{equation}
From this it follows easily that the torsion coefficients $t_i(K)$ are 
given by
\begin{equation}\label{eqn:L-knot-torsion}
\begin{array}{ll}
t_i(K)=n_k-n_{k-1}+\cdots+n_{k-2j}-i &\ \ \ n_{k-2j-1}\le i\le n_{k-2j}\\
t_i(K)=n_k-n_{k-1}+\cdots+n_{k-2j}-n_{k-2j-1} &\ \ \ n_{k-2j-2}\le i\le n_{k-2j-1},\\
\end{array}
\end{equation}
where $j=0,\ldots,k-1$ and $n_{-j}=-n_j$, $n_0=0$. In particular, $t_i(K)$ is
nonincreasing in $i$ for $i\ge 0$. 

The following formula for $d$-invariants 
of surgeries on such a knot is based on results in \cite{Zsurg} and \cite{lenssurg}
(see also \cite[Theorem 1.2]{Qsurg}).  The proof was outlined to us by Peter \ozsvath.

\begin{theorem}
\label{thm:L-space-knot-d}
Let $K\subset S^3$ be an $L$-space knot and let $t_i(K)$ denote its torsion coefficients. 
Then for any $n>0$ the $d$-invariants of the $\pm n$ surgery on $K$ are given by
$$d(K_n,i)=d(U_n,i)-2t_i(K),\ \ \ d(K_{-n},i)=-d(U_n,i)$$
for $|i| \le n/2$, where $U_n$ denotes the $n$ surgery on the unknot $U$. 
\end{theorem}

Before sketching a proof of the theorem we need to explain the notation.
The $d$-invariants are usually associated to \spinc structures on the manifold. The
set of \spinc structures on a three-manifold $Y$ is parametrized by $H^2(Y;\zz)\cong H_1(Y;\zz)$. 
In the case of interest we have $H_1(K_{\pm n};\zz)\cong\zz/n\zz$, hence \spinc
structures can be labelled by the elements of $\zz/n\zz$. Such a labelling is assumed 
in the above theorem and described explicitly as follows. Attaching a 2-handle with
framing $n$ to $S^3\times [0,1]$ along $K \subset S^3\times \{1\}$ gives a cobordism
$W_n$ from $S^3$ to $K_n$. Note that $H_2(W_n;\zz)$ is generated by the homology class
of the core of the 2-handle attached to a Seifert surface for $K$; denote this generator
by $f_n$. A \spinc structure $\spinct$ on $K_n$ is labelled by $i$ if it admits an extension
$\spincs$ to $W_n$ satisfying
$$\langle c_1(\spincs),f_n\rangle -n\equiv 2i \ \pmod{2n}.$$

If $K=U$ is the unknot, the surgery is the lens space which bounds the disk bundle over $S^2$ with Euler number $n$.   
Using either \cite[Proposition 4.8]{os4} or \cite[Corollary 1.5]{os6}, the $d$-invariants 
of $U_n$ for $n>0$ and $|i|<n$ are given by
$$d(U_n,i)=\frac{(n-2|i|)^2}{4n}-\frac14\, .$$

\proof[Proof of Theorem \ref{thm:L-space-knot-d}.]
According to \cite[Theorem 4.1]{Zsurg} the Heegaard Floer homology groups $HF^+(K_n,i)$
for $n\ne 0$ can be computed from the knot Floer homology of $K$. The knot complex
$C=CFK^\infty(S^3,K)$ is a $\zz^2$-filtered chain complex, which is a finitely generated
free module over $\calt:=\zz[U,U^{-1}]$. We write $(i,j)$ for the components of bidegree.
Here $U$ denotes a formal variable in degree $-2$ that decreases the bidegree by $(1,1)$.
The homology of the complex $C$ is $HF^\infty(S^3)=\calt$, the homology of the quotient
complex $B^+:=C\{i\ge 0\}$ is $HF^+(S^3)=\zz[U^{-1}]=:\calt^+$ and the homology of
$C\{i=0\}$ is $\HFhat(S^3)=\zz$. The complexes $C\{j\ge 0\}$ and $B^+$ are quasi-isomorphic
and we fix a chain homotopy equivalence from $C\{j\ge 0\}$ to $B^+$.

We now recall the description of $HF^+(K_{\pm n})$ ($n>0$) in the case of an $L$-space knot $K$. 
For $s\in\zz$, let $A_s^+$ denote the quotient complex
$C\{i\ge 0 \text{\ or\ } j\ge s\}$. 
Let $v_s^+:A_s^+\to B^+$ denote the projection and $h_s^+:A_s^+\to B^+$ the
chain map defined by first projecting to $C\{j\ge s\}$, then applying $U^s$ to identify with
$C\{j\ge 0\}$ and finally applying the chain homotopy equivalence to $B^+$. For any $\sigma\in 
\{0,1,\ldots,n-1\}$ let $\aa^+_\sigma=\oplus_{s\in\sigma+n\zz}A_s^+$ and 
$\bb^+_\sigma=\oplus_{s\in\sigma+n\zz}B^+_s$, where $B_s^+=B^+$ for all $s\in\zz$.
Let $D_{\pm n}^+:\aa^+_\sigma\to\bb^+_\sigma$ be a homomorphism that on $A_s^+$ acts by 
$D_{\pm n}^+(a_s)=(v_s^+(a_s),h_s^+(a_s)) \in B_s^+\oplus B_{s\pm n}^+$. Then $HF^+(K_{\pm n},\sigma)$
is isomorphic to the direct sum of the kernel and the cokernel of the map that $D_{\pm n}^+$
induces between the homologies of $\aa^+_\sigma$ and $\bb^+_\sigma$ (which are direct sums of 
the homologies of $A_s^+$ and $B_s^+$). Moreover, this isomorphism is a homogeneous map of
degree $\pm d(U_n,\sigma)$, where $\aa^+_\sigma$ is graded so that $D_{\pm n}^+$ has degree $-1$. When
computing for $K_n$, the grading on $B^+_s$, where $s=\sigma+nk$, is such that 
$U^0\in B^+_s$ has grading $2k\sigma+nk(k-1)-1$. In case of $-n$-surgery 
$U^0\in B^+_{\sigma+kn}$ has grading $-2(k+1)\sigma-nk(k+1)$. (See \cite{Zsurg} for more details.)

\begin{figure}[htbp]
  \begin{center}
\ifpic
    \leavevmode
    \[ \xy/r1.8pc/:
    (-4,4)*{}; (4,4)*{} **\dir{.}; 
    (-4,3)*{}; (4,3)*{} **\dir{.}; 
    (-4,2)*{}; (4,2)*{} **\dir{.}; 
    (-4,1)*{}; (4,1)*{} **\dir{.};
    {\ar@{.>} (-5,0)*{}; (5,0)*{}}; 
    (-4,-1)*{}; (4,-1)*{} **\dir{.}; 
    (-4,-2)*{}; (4,-2)*{} **\dir{.}; 
    (-4,-3)*{}; (4,-3)*{} **\dir{.}; 
    (-4,-4)*{}; (4,-4)*{} **\dir{.};
    (-4,-4)*{}; (-4,4)*{} **\dir{.}; 
    (-3,-4)*{}; (-3,4)*{} **\dir{.}; 
    (-2,-4)*{}; (-2,4)*{} **\dir{.}; 
    (-1,-4)*{}; (-1,4)*{} **\dir{.};
    {\ar @{.>}(0,-5)*{}; (0,5)*{}}; 
    (1,-4)*{}; (1,4)*{} **\dir{.}; 
    (2,-4)*{}; (2,4)*{} **\dir{.}; 
    (3,-4)*{}; (3,4)*{} **\dir{.}; 
    (4,-4)*{}; (4,4)*{} **\dir{.};
    {\ar(1,4)*{\bullet}; (0,4)*{\bullet}}; 
    {\ar(4,1)*{\bullet}; (4,0)*{\bullet}}; 
    {\ar(1,4)*{\bullet}; (1,2)*{\bullet}}; 
    {\ar(4,1)*{\bullet}; (2,1)*{\bullet}}; 
    {\ar(2,2)*{\bullet}; (2,1)*{\bullet}}; 
    {\ar(2,2)*{\bullet}; (1,2)*{\bullet}}; 
    {\ar(0,3)*{\bullet}; (-1,3)*{\bullet}}; 
    {\ar(3,0)*{\bullet}; (3,-1)*{\bullet}}; 
    {\ar(0,3)*{\bullet}; (0,1)*{\bullet}}; 
    {\ar(3,0)*{\bullet}; (1,0)*{\bullet}}; 
    {\ar(1,1)*{\bullet}; (1,0)*{\bullet}}; 
    {\ar(1,1)*{\bullet}; (0,1)*{\bullet}}; 
    {\ar(-1,2)*{\bullet}; (-2,2)*{\bullet}}; 
    {\ar(2,-1)*{\bullet}; (2,-2)*{\bullet}}; 
    {\ar(-1,2)*{\bullet}; (-1,0)*{\bullet}}; 
    {\ar(2,-1)*{\bullet}; (0,-1)*{\bullet}}; 
    {\ar(0,0)*{\bullet}; (0,-1)*{\bullet}}; 
    {\ar(0,0)*{\bullet}; (-1,0)*{\bullet}}; 
    {\ar(-2,1)*{\bullet}; (-3,1)*{\bullet}}; 
    {\ar(1,-2)*{\bullet}; (1,-3)*{\bullet}}; 
    {\ar(-2,1)*{\bullet}; (-2,-1)*{\bullet}}; 
    {\ar(1,-2)*{\bullet}; (-1,-2)*{\bullet}}; 
    {\ar(-1,-1)*{\bullet}; (-1,-2)*{\bullet}}; 
    {\ar(-1,-1)*{\bullet}; (-2,-1)*{\bullet}}; 
    {\ar(-3,0)*{\bullet}; (-4,0)*{\bullet}}; 
    {\ar(0,-3)*{\bullet}; (0,-4)*{\bullet}}; 
    {\ar(-3,0)*{\bullet}; (-3,-2)*{\bullet}}; 
    {\ar(0,-3)*{\bullet}; (-2,-3)*{\bullet}}; 
    {\ar(-2,-2)*{\bullet}; (-2,-3)*{\bullet}}; 
    {\ar(-2,-2)*{\bullet}; (-3,-2)*{\bullet}}; 
    {\ar(2,5)*{}; (2,3)*{\bullet}}; 
    {\ar(5,2)*{}; (3,2)*{\bullet}}; 
    {\ar(3,3)*{\bullet}; (3,2)*{\bullet}}; 
    {\ar(3,3)*{\bullet}; (2,3)*{\bullet}}; 
    {\ar(3,6)*{}; (3,4)*{\bullet}}; 
    {\ar(6,3)*{}; (4,3)*{\bullet}}; 
    {\ar(4,4)*{\bullet}; (4,3)*{\bullet}}; 
    {\ar(4,4)*{\bullet}; (3,4)*{\bullet}}; 
    {\ar(-4,-1)*{\bullet}; (-5,-1)*{}}; 
    {\ar(-1,-4)*{\bullet}; (-1,-5)*{}}; 
    {\ar(-4,-1)*{\bullet}; (-4,-3)*{\bullet}}; 
    {\ar(-1,-4)*{\bullet}; (-3,-4)*{\bullet}}; 
    {\ar(-3,-3)*{\bullet}; (-3,-4)*{\bullet}}; 
    {\ar(-3,-3)*{\bullet}; (-4,-3)*{\bullet}}; 
    {\ar(-4,-4)*{\bullet}; (-4,-5)*{}}; 
    {\ar(-4,-4)*{\bullet}; (-5,-4)*{}}; 
    (5.5,0)*{i};
    (0,5.5)*{j};
    (5,5)*{.};(6,6)*{.};(7,7)*{.};
    (-5,-5)*{.};(-6,-6)*{.};(-7,-7)*{.};
\endxy \]
\else \vskip 5cm \fi
   \begin{narrow}{0.3in}{0.3in}
    \caption{{\bf The knot Floer complex $CFK^\infty(T_{3,5})$.}  Each bullet represents
    a $\zz$, and each arrow is an isomorphism.  The groups and differentials on the axes are determined by the Alexander polynomial
    $$\Delta_{T_{3,5}}(T)=-1+(T+T^{-1})-(T^3+T^{-3})+(T^4+T^{-4}),$$
    and these in turn determine the entire complex.}
    \label{fig:CFK}
    \end{narrow}
  \end{center}
\end{figure}

Suppose now that $K$ is an $L$-space knot with Alexander polynomial as in 
(\ref{eqn:L-knot-alex}). Define $\delta_k:=0$ and
$$\delta_l:=\left\{\begin{array}{ll}
\delta_{l+1}-2(n_{l+1}-n_l)+1 & \text{\ \ \ if $k-l$ is odd}\\
\delta_{l+1}-1 & \text{\ \ \ if $k-l$ is even}
\end{array}\right.
$$
for $l=k-1,k-2,\ldots,-k$, where as above $n_{-l}=-n_l$. Then by \cite[Theorem 1.2]{lenssurg}
$C\{i=0\}$ is (up to quasi-isomorphism) equal to the free abelian group with one generator $x_l$
in bidegree $(0,n_l)$ for $l=-k,\ldots,k$ and the grading of $x_l$ is $\delta_l$. To determine
the differentials note that the homology of $C\{i=0\}$ is $\zz$ in grading 0, so generated by 
(the homology class of) $x_k$. It follows that the differential on $C\{i=0\}$ is a collection
of isomorphisms $C_{0,n_{k-2l+1}}\to C_{0,n_{k-2l}}$ for $l=1,\ldots,k$. Similarly we see
that the differential on $C\{j=0\}$ is given by a collection of isomorphisms 
$C_{n_{k-2l+1},0}\to C_{n_{k-2l},0}$ for $l=1,\ldots,k$. This, together with U-equivariance,
determines the complex $C$.  (For an example, see Figure \ref{fig:CFK}.)

Suppose that the homology $H$ of a quotient complex of $C$ is isomorphic
to $\calt^+=\zz[U^{-1}]$.  We say $H$ starts at $(i,j)$ if the element $U^0$ 
has a representative of bidegree 
$(i,j)$. With this notation $H_*(B_s)$ starts at $(0,n_k)$ and
$H_*(C\{j\ge s\})$ starts at $(s,s+n_k)$. It remains to consider $H_*(A_s^+)$. Note that
these groups are also isomorphic to $\calt^+$. For $s\ge n_k$ the homology of  $A_s^+$ starts at $(0,n_k)$,
so $v_s^+=id$ and $h_s^+=U^s$. For $n_{k-1}\le s < n_k$, $H_*(A_s^+)$ starts at $(s-n_k,s)$, so
$v_s^+=U^{n_k-s}$ and $h_s^+=U^{n_k}$. For $n_{k-2}\le s<n_{k-1}$ it starts at
$(n_{k-1}-n_k,n_{k-1})$, hence $v_s^+=U^{n_k-n_{k-1}}$ and $h_s^+=U^{n_k-n_{k-1}+s}$.
It is now easy to observe that for $s\ge 0$ the homology $H_*(A_s^+)$ starts at
$(t_s,t_s+n_k)$ and thus $v_s^+=U^{t_s}$ and $h_s^+=U^{t_s+s}$, where $t_s$ is a torsion
coefficient of $K$ (compare with equation (\ref{eqn:L-knot-torsion})). Similarly for
$s<0$ we obtain $v_s^+=U^{t_s-s}$ and $h_s^+=U^{t_s}$. 

Note that since the \spinc structures corresponding to $i$ and $-i$ are conjugate (so their $d$ invariants are equal) we may restrict to those in the range $\{0,1,\ldots,\lfloor n/2\rfloor\}$.
Consider now $K_n$ ($n>0$) and choose some $\sigma\in \{0,1,\ldots,\lfloor n/2\rfloor\}$. Writing
$D_n^+:\aa^+_\sigma\to\bb^+_\sigma$ in components we obtain for $l>0$:
$$\begin{array}{l}
b_{\sigma+ln}=U^{t_{\sigma+ln}}a_{\sigma+ln}+U^{t_{\sigma+(l-1)n}+\sigma+(l-1)n}a_{\sigma+(l-1)n}\\
\qquad  \qquad \qquad \implies 
U^{t_{\sigma+ln}}a_{\sigma+ln}=b_{\sigma+ln}-U^{t_{\sigma+(l-1)n}+\sigma+(l-1)n}a_{\sigma+(l-1)n}\\
b_\sigma=U^{t_\sigma}a_\sigma+U^{t_{n-\sigma}}a_{\sigma-n}\\
\qquad  \qquad \qquad \implies 
U^{t_{n-\sigma}}a_{\sigma-n}=b_\sigma-U^{t_\sigma}a_\sigma\\ 
b_{\sigma-ln}=U^{t_{ln-\sigma}+ln-\sigma}a_{\sigma-ln}+U^{t_{(l+1)n-\sigma}}a_{\sigma-(l+1)n}\\
\qquad  \qquad \qquad \implies 
U^{t_{(l+1)n-\sigma}}a_{\sigma-(l+1)n}=b_{\sigma-ln}-U^{t_{ln-\sigma}+ln-\sigma}a_{\sigma-ln}
\end{array}
$$
from which it is easily seen that $D_n^+$ is surjective and its kernel contains one $\calt^+$ summand, 
isomorphic to $H_*(A_{\sigma})$. Since $D_n^+$ shifts grading by  $-1$, $U^0 \in H_*(B_\sigma)$ has grading $-1$, and the component of $D_n^+$ from $H_*(A_\sigma)$ to $H_*(B_\sigma)$ 
is equal to $U^{t_\sigma}$, it follows 
that $U^0\in \calt^+\subset \ker D_n^+$ has grading $-2t_\sigma$.
The formula for $d(K_n,\sigma)$ now follows using the degree shift between $HF^+(K_n,\sigma)$
and $\ker D_n^+$.

Finally consider $K_{-n}$ ($n>0$) and choose $\sigma\in \{0,1,\ldots,\lfloor n/2\rfloor\}$. Writing
$D_{-n}^+$ in components yields:
$$\begin{array}{l}
b_{\sigma+ln}=U^{t_{\sigma+ln}}a_{\sigma+ln}+U^{t_{\sigma+(l+1)n}+\sigma+(l+1)n}a_{\sigma+(l+1)n} \hfill  (l\ge 0)\\
\qquad  \qquad \qquad \implies 
U^{t_{\sigma+ln}}a_{\sigma+ln}=b_{\sigma+ln}-U^{t_{\sigma+(l+1)n}+\sigma+(l+1)n}a_{\sigma+(l+1)n}\\
b_{\sigma-n}=U^{t_{\sigma}+\sigma}a_\sigma+U^{t_{n-\sigma}+n-\sigma}a_{\sigma-n}\\
\qquad  \qquad \qquad \implies 
U^{t_{\sigma}+\sigma}a_\sigma=b_{\sigma-n}-U^{t_{n-\sigma}+n-\sigma}a_{\sigma-n}\\
b_{\sigma-ln}=U^{t_{ln-\sigma}+ln-\sigma}a_{\sigma-ln}+U^{t_{(l-1)n-\sigma}}a_{\sigma-(l-1)n}  \hfill (l\ge 2)\\
\qquad  \qquad \qquad \implies 
U^{t_{(l-1)n-\sigma}}a_{\sigma-(l-1)n}=b_{\sigma-ln}-U^{t_{ln-\sigma}+ln-\sigma}a_{\sigma-ln}\, .
\end{array}
$$
The top and bottom equations determine $a_s$ (or more precisely $U^N a_s$ for $N>>0$) for all $s\in\sigma+n\zz$, so the middle equation cannot be fulfilled 
in general. It follows that $D_{-n}^+$ has cokernel isomorphic to $H_*(B_{\sigma-n})=\calt^+$; the grading of $U^0$
in this group is $0$.  The formula for $d(K_{-n},\sigma)$ again follows from the degree shift.
\endproof

Combining Theorems \ref{thm:YboundsX} and \ref{thm:L-space-knot-d} yields 

\begin{theorem}
\label{thm:L-space-knot-surgery}
Let $n>0$ and let $K$ be an $L$-space knot whose torsion coefficients satisfy
$$t_i(K)>\left\{\begin{array}{ll}
\frac{(n-2i)^2+1}{8n} -\frac14 & \text{if $n$ is odd}\\
\frac{(n-2i)^2}{8n} -\frac14 & \text{if $n$ is even}
\end{array}\right.
$$
for $0\le i\le n/2$.
Then for $0<m\le n$, $K_m$ cannot bound a negative-definite four-manifold with no torsion in $H_1$. 
\end{theorem}

\proof
The formulas follow from the above-mentioned Theorems. To see that obstruction for $n$-surgery to bound 
a negative-definite manifold implies obstruction for all $m$-surgeries with $0<m\le n$ 
observe that for fixed $i$ the right-hand side of the inequality is an increasing function of $n$
for $n\ge 2i$.
\endproof

Note that the surgery coefficient $m$ in Theorem \ref{thm:L-space-knot-surgery} is an integer.
In \cite{bsinprep} we will consider in more detail the question of which surgeries 
(including Dehn surgeries) on
a knot $K$ can be given as the boundary of a negative-definite four-manifold.  In particular
we will show that Theorem \ref{thm:L-space-knot-surgery} holds as stated with $m\in\qq$.



\section{Example: Surgeries on torus knots}
\label{sec:torus-knots}

In this section we consider torus knots $T_{p,q}$; we assume $2\le p<q$.  Right-handed 
torus knots are $L$-space knots since for example $pq-1$ surgery on $T_{p,q}$ yields
a lens space \cite{moser}. 
Let $N=(p-1)(q-1)/2$ denote the degree of the Alexander polynomial of $T_{p,q}$.  The following proposition gives a simple function which approximates the torsion coefficients of a torus knot.

\begin{proposition}
\label{prop:torus-torsion}
The torsion coefficients of $T_{p,q}$ are given by
$$t_i=\#\{(a,b)\in\zz_{\ge0}^2\, |\, ap+bq<N-i\}$$
and they satisfy $t_i\ge g(N-i)$ for $0\le i\le N$, where
$x\mapsto g(x)$ is a piecewise linear continuous function, which equals $0$ for $x\le0$ and
whose slope on the interval $[(k-1)q,kq]$ is $k/p$.
\end{proposition}

\proof
The (unsymmetrised) Alexander polynomial of $K=T_{p,q}$ is
$$\tilde\Delta_K(T)=\frac{(1-T^{pq})(1-T)}{(1-T^p)(1-T^q)}\, ,$$
which is a polynomial of degree $2N$. Writing $\tilde\Delta_K$ as a formal power series in $T$
we obtain
\begin{eqnarray*}
\tilde\Delta_K(T)&=&(1-T^{pq})(1-T)\sum_{a\ge0}T^{ap}\sum_{b\ge0}T^{bq}\hfill\\
&=&(1-T)\sum_{a\ge0}T^{ap}\sum_{b\ge0}T^{bq} - T^{pq}(1-T)\sum_{a\ge0}T^{ap}\sum_{b\ge0}T^{bq}\, .
\end{eqnarray*}
Clearly only the terms in the first product of the last line contribute to the nonzero coefficients
of $\Delta_K$. A power $T^k$ appears with coefficient $+1$ in this term whenever $k=ap+bq$ for some
$a,b\in\zz_{\ge0}$, and with coefficient $-1$ whenever $k-1=ap+bq$ for some $a,b\in\zz_{\ge0}$.
Since the coefficient $a_j$ in $\Delta_K(T)$ is the coefficient of $T^{N-j}$ in 
$\tilde\Delta_K(T)$, we obtain $a_j=m_j-m_{j+1}$, where
$$m_j:=\#\{(a,b)\in\zz_{\ge0}^2\, |\, ap+bq=N-j\}\, .$$
Then for $i\ge 0$
$$t_i=\sum_{j>0}ja_{i+j}=\sum_{j>0}m_{i+j}=\#\{(a,b)\in\zz_{\ge0}^2\, |\, ap+bq<N-i\}\, .$$
In what follows it is convenient to replace $t_i$ by 
$$s(i):=t_{N-i}=\#\{(a,b)\in\zz_{\ge0}^2\, |\, ap+bq<i\}.$$
Define 
$$
\bar g(i):=\left\{\begin{array}{ll}
i/p& \text{if } i\ge 0\cr
0 & \text{if } i\le 0
\end{array}\right. \, 
$$
and $\bar s(i):=\ceil{g(i)}$; clearly $\bar s(i)\ge \bar g(i)$ for every $i$. 
Separating the set appearing in the definition of $s(i)$ into subsets with fixed value of $b$ we
obtain
$$s(i)=\sum_{b\ge0}\#\{a\in\zz_{\ge0}\, |\, ap<i-bq\}=
\sum_{b\ge0}\bar s(i-bq) \ge \sum_{b\ge0} \bar g(i-bq)=:g(i)\, .$$
\endproof

The following corollary describes the range of surgeries on $T_{p,q}$ that cannot bound 
negative-definite manifolds according to Theorem \ref{thm:L-space-knot-surgery} (and using Proposition \ref{prop:torus-torsion}). To obtain
the result of Proposition \ref{prop:torus-neg-def} note that all the lower bounds in the
corollary allow for $m=2$.  Note that Lisca and Stipsicz have shown that Dehn surgery on $T_{2,2n+1}$
with positive framing $r$ bounds a negative-definite manifold (possibly with torsion in $H_1$) if and
only if $r\ge 4n$ \cite{ls}.

\begin{corollary}
\label{cor:torus-neg-def}
Let $2\le p<q$ and $N=(p-1)(q-1)/2$.  
\begin{itemize}
\item
If $p$ is even and $m$ is less than the minimum of
$$
\begin{array}{l}
1+\sqrt{4N}\cr
2+\frac12\sqrt{\alpha(\alpha+4\beta)-4}-\frac12\alpha\cr
q-p+3
\end{array}
$$
where
$\alpha=q(p-2)+2$ and $\beta=q-p+1$ 
\item
or if $p$ is odd and $m$ is less than
the minimum of 
$$
\begin{array}{l}
1+\sqrt{4N}\cr
2+\frac12\sqrt{\alpha(\alpha+4\beta)-4}-\frac12\alpha-\frac{q(p-3)}p\cr
q-p+5-\frac{q+2}p
\end{array}
$$
where
$\alpha=q(p-4)+2+3q/p$ and $\beta=2q-p+1$ 
\item
and $1\le n\le 2N+m$, 
\end{itemize}
then $+n$-surgery on $T_{p,q}$ cannot bound a negative definite four-manifold
with no torsion in $H_1$.
\end{corollary}

\proof
Write $n=2N+m$; we may assume $n<pq-1$. It suffices to show that
\begin{equation}\label{eqn:hg}
h(x):=\frac{(m+2x)^2+1}{8n}-\frac14<g(x)
\end{equation}
for $-m/2\le x\le N$, where $g$ is the function appearing in Proposition \ref{prop:torus-torsion}.

Since $h$ is convex and $g$ is piecewise linear, it is enough to check the inequality 
for $x=kq$ with $k=0,1,\ldots,\floor{N/q}$ and for $x=N$. Consider first $x=kq$. 
From the definition of $g$ we obtain $g(kq)=qk(k+1)/2p$. Substituting this into 
(\ref{eqn:hg}) we obtain
$$4k^2(q^2-n\frac qp)+4k(qm-n\frac qp)+m^2+1-2n<0\, .$$
Since $q-n/p>1/p$, it suffices to consider only $k=0$ and $k=\floor{N/q}$. For 
$k=0$ the last inequality yields $m<1+\sqrt{4N}$.

Assume first $p$ is even; then $\floor{N/q}=p/2-1$ and substituting this into 
(\ref{eqn:hg}) gives
$$(m+q(p-2))^2+1<n(2+q(p-2)).$$
Writing $\mu=m+q(p-2)$, $\alpha=q(p-2)+2$ and $\beta=q-p+1$ the last inequality
becomes
$$\mu^2-\alpha\mu+1-\alpha\beta<0\, ,$$
which implies 
$$\mu<\frac{\alpha}2+\frac12\sqrt{\alpha(\alpha+4\beta)-4}\, ;$$
this is equivalent to the second condition on $m$ in the statement of the corollary.

Since the slope of $g$ on the interval from $(p/2-1)q$ to $N$ is $1/2$, 
we get
$g(N)=g((p/2-1)q)+(N-(p/2-1)q)/2=(pq-2p+2)/8$. Subsituting this into (\ref{eqn:hg}) gives
$$n^2-n(pq-2p+4)+1<0\, ,$$
which holds if $n<pq-2p+4$ or $m<q-p+3$.

Assume now $p$ is odd; then $\floor{N/q}=(p-3)/2$ and substituting this into 
(\ref{eqn:hg}) gives
$$(m+q(p-3))^2+1<n(2+\frac qp (p-1)(p-3)).$$
Writing $\mu=m+q(p-3)$, $\alpha=q(p-4)+2+3q/p$ and $\beta=2q-p+1$ the last inequality
becomes
$$\mu^2-\alpha\mu+1-\alpha\beta<0\, ,$$
from which the second condition on $m$ in the statement of the corollary follows.

Now the slope of $g$ on the interval from $q(p-3)/2$ to $N$ is $(p-1)/2p$ and
$g(N)=(pq-2p+4-(q+2)/p)/8$. Subsituting this into (\ref{eqn:hg}) gives
$$n^2-n(pq-2p+6-\frac{q+2}p)+1<0\, ,$$
which holds if $n<pq-2p+6-\frac{q+2}p$ or $m<q-p+5-\frac{q+2}p$.
\endproof



\begin{thebibliography}{99}
\bibitem{elkies} N.~D.~Elkies.  \textsl{A characterization of the $ Z\sp n$ lattice},  
Math.~Res.~Lett.  {\bf 2}  1995,  321--326.
\bibitem{f} K.~A.~Fr{\o}yshov.  \textsl{The Seiberg-Witten equations and 
four-manifolds with boundary},  Math. Res. Lett.  {\bf 3}  1996,  373--390.
\bibitem{kk} A.~Kawauchi \& S.~Kojima. \textsl{Algebraic classification of
linking pairings on 3-manifolds},
Math.~Ann.~ {\bf 253} 1980, 29--42.
\bibitem{kp} M.~Kneser \& D.~Puppe. \textsl{Quadratische Formen und
Verschlingungsinvarienten von Knoten},
Math.~Zeitschr.~ {\bf 58} 1953, 376--384.
\bibitem{ls} P.~Lisca \& A.~I.~Stipsicz. \textsl{\ozsvath-\szabo\ invariants
and tight contact structures,~I}, Geometry and Topology {\bf 8} 2004, 925--945.
\bibitem{moser} L.~Moser. \textsl{Elementary surgery along a torus knot},
    Pacific Journal of Math. {\bf 38} 1971, 737--745.
\bibitem{bs} B.~Owens \& S.~Strle.  \textsl{Rational homology
spheres and the four-ball genus of knots},  Advances in Mathematics {\bf 200} 2006, 196--216.
\bibitem{bsinprep}B.~Owens \& S.~Strle.  \textsl{Dehn surgeries and negative-definite four-manifolds}, in preparation.
\bibitem{os4} P.~\ozsvath~ \& Z.~\szabo. \textsl{Absolutely graded Floer
    homologies and intersection forms for four-manifolds with boundary},
    Advances in Mathematics {\bf 173} 2003, 179--261.
\bibitem{os6} P.~\ozsvath~ \& Z.~\szabo. \textsl{On the Floer
    homology of plumbed three-manifolds},
    Geometry and Topology {\bf 7} 2003, 225--254.
\bibitem{lenssurg} P.~\ozsvath~ \& Z.~\szabo. \textsl{On knot Floer
    homology and lens space surgeries},  Topology {\bf 44}  2005, 1281--1300.
\bibitem{Zsurg} P.~\ozsvath~ \& Z.~\szabo. \textsl{Knot Floer
    homology and integer surgeries}, math.GT/0410300.
\bibitem{Qsurg} P.~\ozsvath~ \& Z.~\szabo. \textsl{Knot Floer
    homology and rational surgeries}, math.GT/0504404.
\bibitem{serre} J.-P.~Serre. \textsl{A course in arithmetic}, Springer, 1973.
\bibitem{taylor} L.~R.~Taylor. \textsl{Relative Rochlin invariants}, 
Topology Appl. {\bf 18} 1984, 259--280.
\bibitem{wall} C.~T.~C.~Wall. \textsl{Quadratic forms on finite groups, and 
related topics}, Topology {\bf 2} 1963, 281--298.
\end{thebibliography}
\end{document}